\documentclass[10pt]{article}
\usepackage{graphicx}
\usepackage[body={6.0in, 8.2in},left=1.25in,right=1.25in]{geometry}
\usepackage{amsfonts}
\usepackage{amssymb}
\usepackage{amsmath}
\usepackage{multicol}
\usepackage{subfigure}

\begin{document}
\title{ On convergence and stability of a numerical scheme of coupled nonlinear Schr\"{o}dinger equations}
\author{Bartosz Reichel, Sergey Leble}
\maketitle
\begin{abstract}
We consider numerical solution of  Coupled Nonlinear Schr\"{o}dinger
Equation. We prove stability and convergence in the $L_2$ space for
an explicit scheme which estimations is used for implicit scheme and
compare both method. As a test  we compare numerical
 solutions of Manakov system with known analytical solitonic
 solutions and as
  example of general system - evolution of two impulses with
different group velocity (model of pulses interaction in optic
fibers). Last example,  a rectangular pulse evolution, shows
asymptotic behavior typical for Nonlinear Schr\"{o}dinger Equation
asymptotics with the same initial condition.
\end{abstract}

\section{Introduction}
A growing interest to encoded electromagnetic pulses propagation for
long distances and processing is noted
\cite{book:Mollenauer:OpticalSolitons,book:agrawal:NonFibOpt}. A lot
of realistic models are based on the Coupled Nonlinear
Schr\"{o}dinger Equations (CNLSE): they are developed and used in
many nonlinear optics problems such as polarization modes
interaction \cite{book:agrawal:NonFibOpt,Menyuk1987,GarthPask1992}.
For such modes some numerical schemes and codes are elaborated
beginning from a celebrated integrable Ablowitz-Ladik one
\cite{Ablowitz1976a,Ablowitz1976b}.

In papers
\cite{Chang1999,Ismail2004,Ivanauska1999,Sun2004,
Kurtinaitis2004,LWu1996} authors investigate numerical methods for
solving CNLS equation based on finite difference schemes. Such
difference schemes are applied in \cite{Ablowitz2002} to model
vector spatial soliton behavior in a nonlinear waveguide arrays.

There are numerical schemes which are unconditionally unstable,
conditionally stable and unconditionally stable. Generally could be
said that we have two useful schemes: conditionally stable (Euler
type) schemes and  unconditionally stable (Crank-Nicolson type)
schemes; we focus on these types of schemes. Most important thing is
to use conservation laws for CNLS equations while the scheme
constructed. There are a lot of possibilities to define conservation
laws \cite{Ismail2004} but we choose in this paper the standard one
given by \cite{Sun2004,Ivanauska1999}.

Authors of \cite{Ivanauska1999,Kurtinaitis2004,Chang1999} compare
explicit method with implicit one underlying that both of the
methods have good and bad sides and there is no an uniform point of
view. Conditionally stable methods do not need matrix formulation,
but they need smaller space and time steps to assure stability.
Unconditionally stable method need to solve system of equations, but
bigger time step could be used. We would like to compare numerical
solutions for implicit and explicit method under the scope of
stability proof for explicit method. In paper \cite{Ivanauska1999}
authors prove convergence and stability for NLS equation but do not
analyze CNLS equations and how stability and convergence depends on
initial amplitude (energy) of impulses.

 We consider a finite-difference scheme for the CNLSE which extends in a sense the results
concerned linear Schr\"{o}dinger equation \cite{WDai1992}, following
the ideas  successfully applied for the coupled Korteveg -de Vries
(KdV) systems in \cite{Halim2002,Halim2003}.

Most important aim which we would like to achieve is to estimate
stability and convergence parameters as a function of initial
amplitude of pulses (energy) for the CNLS equations on a base of the
norm originated from $L_2$ space. Such result allows us to estimate
time and space steps for the convergence regime inside which we
could consider a solution as stable.

The system of the CNLSE is considered in the form:
\begin{subequations}\label{E:CNLSE_2}
\begin{equation}
i\partial_{t}U+ i\sigma\partial_{x}U+k
\partial_{xx}U+\left[a|U|^{2}+b|V|^{2}\right]U=0,
\end{equation}
\begin{equation}
i\partial_{t}V-i\sigma\partial_{x}V+k
\partial_{xx}V+\left[c|V|^{2}+d|U|^{2}\right]V=0,
\end{equation}
\end{subequations}
where functions $U$ and $V$ are differentiable up to the second
order and the following notations are used
\begin{subequations}
\begin{eqnarray}
U&=&U(x,t),\\
V&=&V(x,t),\\
\partial_{t}V\equiv U_{t}&\equiv&\frac{\partial U(x,t)}{\partial t},\\
\partial_{xx}V\equiv U_{xx}&\equiv&\frac{\partial^{2} U(x,t)}{\partial x^{2}},
\end{eqnarray}
\end{subequations}
 and (for the example of two modes excited in the waveguide \cite{Snyder1978,Leble2005}) the  parameter
 $\sigma$ is defined as
\begin{equation}
\sigma=\frac{1}{2}\left(k'_{01}-k'_{11}\right).
\end{equation}
This parameter describes the difference between   group velocities
of the modes. For simplicity of notation, but without losses of
generality, we take $a,b,c,d,k,\sigma \geq 0$ (in other case we
should put everywhere in estimation of norm absolute value of this
coefficients which could make equations more illegible).

\section{Conservation laws}

Let us  prove first that the solutions of the system
\eqref{E:CNLSE_2} satisfy the following conservation laws
\begin{subequations}
\begin{eqnarray}
\int\limits_{-\infty}^{\infty}|U|^{2}d\tau=const,\\
\int\limits_{-\infty}^{\infty}|V|^{2}d\tau=const.
\end{eqnarray}
\end{subequations}

If one writes the conjugate to the first of equations
\eqref{E:CNLSE_2} (equation for $U$ amplitude)
\begin{equation}\label{E:con_CNLSE_2_1}
-i\partial_{\xi}\overline{U}+i\delta\partial_{\tau}\overline{U}+k
\partial_{\tau\tau}\overline{U}+\left[a|U|^{2}+b|V|^{2}\right]\overline{U}=0,
\end{equation}
next multiply the result \eqref{E:con_CNLSE_2_1} by $-U$ and the
equation \eqref{E:CNLSE_2} by $\overline{U}$, it yields
\begin{subequations}
\begin{equation}\label{E:cs1}
i\overline{U}\partial_{\xi}U-
i\delta\overline{U}\partial_{\tau}U+k
\overline{U}\partial_{\tau\tau}U+\left[a|U|^{2}+b|V|^{2}\right]\overline{U}U=0,
\end{equation}
\begin{equation}\label{E:cs2}
iU\partial_{\xi}\overline{U}-i\delta
U\partial_{\tau}\overline{U}-k
U\partial_{\tau\tau}\overline{U}-\left[a|U|^{2}+b|V|^{2}\right]U\overline{U}=0,
\end{equation}
\end{subequations}
with obvious relations
\begin{equation}
 U\overline{U}=\overline{U}U=|U|^{2}.
\end{equation}
As a next step,  adding the equation \eqref{E:cs1} to \eqref{E:cs2}
\begin{equation}
i\partial_{\xi}U\overline{U}-i\delta
\partial_{\tau}U\overline{U}+k\partial_{\tau}\left[
\overline{U}\partial_{\tau}U-U\partial_{\tau}\overline{U}\right]=0,
\end{equation}
and integrating the result, one arrives at:
\begin{equation}
i\partial_{\xi}\int\limits_{-\infty}^{\infty}(U\overline{U})d\tau-i\delta
\int\limits_{-\infty}^{\infty}\partial_{\tau}(U\overline{U})d\tau+k\int\limits_{-\infty}^{\infty}\partial_{\tau}\left[
\overline{U}\partial_{\tau}U-U\partial_{\tau}\overline{U}\right]d\tau=
\end{equation}
\begin{equation}
i\partial_{\xi}\int\limits_{-\infty}^{\infty}U\overline{U}d\tau-i\delta
|U|^{2}\Big|^{+\infty}_{-\infty}+k\left[
\overline{U}\partial_{\tau}U-U\partial_{\tau}\overline{U}\right]\Big|^{+\infty}_{-\infty}=0,
\end{equation}

Suppose the boundary conditions at both infinities
\begin{subequations}
\begin{eqnarray}
\lim_{\tau \to \pm \infty} U&=&0,\\
\end{eqnarray}
\end{subequations}
are imposed,
\begin{equation}
i\partial_{\xi}\int\limits_{-\infty}^{\infty}|U|^{2}d\tau=0,
\end{equation}
we obtain the first conservation law in the form
\begin{equation}
\int\limits_{-\infty}^{\infty}|U|^{2}d\tau=const.
\end{equation}
For $V$ we built the conservation law in the same way
\begin{equation}
\int\limits_{-\infty}^{\infty}|V|^{2}d\tau=const.
\end{equation}
Such conservation laws give us a possibility to define the basic
norm in $L_{2}$ space as
\begin{equation}
{\|W\|}^{2}=\int\limits_{-\infty}^{\infty}(|U|^{2}+|V|^{2})d\tau=const.
\end{equation}
We would denote other norms by the same notation if it will not lead
to a confusion.

\section{Discretization of the CNLSE system}

\subsection{Explicit scheme}
Choose the explicit scheme in a simple form
\cite{book:Hoffman:numerical_eng} with a second order discretization
with respect to time and a third order one to space:
\begin{align}\label{LAB:EXPLICITS}
i\frac{U_{i}^{j+1}-U_{i}^{j}}{
\tau}+i\sigma\frac{U_{i+1}^{j}-U_{i-1}^{j}}{2h}+
k\frac{U_{i+1}^{j}-2U_{i}^{j}+U_{i-1}^{j}}{h^{2}}+
\left(a|U_{i}^{j}|^{2}+b|V_{i}^{j}|^{2}\right)U_{i}^{j}=0.
\end{align}

Let us derive conservation law for discrete CNLSE system (where
$\overline{U}$ is complex conjugate of $U$). Multiplying equation
\eqref{LAB:EXPLICITS} by
$(\overline{U}^{j+1}_{i}+\overline{U}^{j}_{i})$  and complex
conjugating of this equation by $(U^{j+1}_{i}+U^{j}_{i})$. If we
apply zero boundary conditions, we have
\begin{subequations}
\begin{eqnarray}\label{E:conserw_law}
\sum_{i=1}^{N}|U_{i}^{j+1}|^{2}&=&\sum_{i=1}^{N}|U_{i}^{j}|^{2},\\
\sum_{i=1}^{N}|V_{i}^{j+1}|^{2}&=&\sum_{i=1}^{N}|V_{i}^{j}|^{2}.
\end{eqnarray}
\end{subequations}

\subsection{Implicit scheme}
This  six point scheme \cite{Sun2004}
\begin{align}
i\frac{U_{l}^{n+1}-U_{l}^{n}}{
\tau}&+i\sigma\frac{U_{l+1}^{n+1}+U_{l+1}^{n}-U_{l-1}^{n+1}-U_{l-1}^{n}}{4h}+
k\frac{U_{l+1}^{n+1}+U_{l+1}^{n}-2U_{l}^{n+1}-2U_{l}^{n}+U_{l-1}^{n+1}+U_{l-1}^{n}}{2h^{2}}+\\
&+\left(a|U_{l}^{n+1/2}|^{2}+b|V_{l}^{n+1/2}|^{2}\right)\frac{U_{l+1}^{n}+U_{l}^{n}}{2}=0,
\end{align}
is exactly Crank-Nicolson one \cite{book:Hoffman:numerical_eng}.
This implicit scheme bases on the same finite difference as explicit
scheme, we expect that implicit scheme converges to exactly solution
quicker than explicit scheme.

Elements in a node $1/2$ are calculated by iterations.

\section{Stability}

We can separate real and imaginary part of $U$ and $V$ by put
\begin{subequations}
\begin{eqnarray}
U = \xi + i\eta,\\
V = \alpha+i\beta.
\end{eqnarray}
\end{subequations}
Substituting this into \eqref{E:CNLSE_2} yields in four equation
with real amplitudes
\begin{subequations}
\begin{eqnarray}
 \xi_{t}+\sigma\xi_{t}+k\eta_{xx}+\left[a(\xi^{2}+\eta^{2})+b(\alpha^{2}+\beta^{2})\right]\eta=0,\\
 -\eta_{t}-\sigma\eta_{t}+k\xi_{xx}+\left[a(\xi^{2}+\eta^{2})+b(\alpha^{2}+\beta^{2})\right]\xi=0,\\
 \alpha_{t}-\sigma\alpha_{t}+k\beta_{xx}+\left[c(\alpha^{2}+\beta^{2})+d(\xi^{2}+\eta^{2})\right]\beta=0,\\
 -\beta_{t}+\sigma\beta_{t}+k\alpha_{xx}+\left[c(\alpha^{2}+\beta^{2})+d(\xi^{2}+\eta^{2})\right]\alpha=0.
\end{eqnarray}
\end{subequations}

If we apply explicit scheme \eqref{LAB:EXPLICITS} to this
equations, we can build time evolution matrix as
\begin{equation}
\mathbf{T}^{j+1} = \left(%
\begin{array}{cccc}
  T_{11}^{j+1} & T_{12}^{j+1} & 0 & 0 \\
  T_{21}^{j+1} & T_{22}^{j+1} & 0 & 0 \\
  0 & 0 & T_{33}^{j+1} & T_{34}^{j+1} \\
  0 & 0 & T_{43}^{j+1} & T_{44}^{j+1} \\
\end{array}%
\right),
\end{equation}
where
\begin{subequations}
\begin{eqnarray}
  T_{11}&=&T_{22} = \delta_{i,r}-\frac{\tau\sigma}{2h}\left(\delta_{i+1,r}-\delta_{i-1,r}\right), \\
  T_{12}&=&-T_{21} = -\frac{\tau k}{h^{2}}\left(\delta_{i+1,r}-2\delta_{i,r}-\delta_{i-1,r}\right)-\tau \delta_{i,r}\left\{a[(\xi^j)^{2}+(\eta^j)^{2}]+b[(\alpha^j)^{2}+(\beta^j)^{2}]\right\}, \\
  T_{33}&=&T_{44} = \delta_{i,r}+\frac{\tau\sigma}{2h}\left(\delta_{i+1,r}-\delta_{i-1,r}\right), \\
  T_{34}&=&-T_{43} = -\frac{\tau k}{h^{2}}\left(\delta_{i+1,r}-2\delta_{i,r}-\delta_{i-1,r}\right)-\tau
  \delta_{i,r}\left\{c[(\alpha^j)^{2}+(\beta^j)^{2}]+d[(\xi^j)^{2}+(\eta^j)^{2}]\right\},
\end{eqnarray}
\end{subequations}
that acts in the vector space $R_w$ of the columns:
\begin{equation}
\mathbf{W}^{j}=\left(%
\begin{array}{c}
  \xi^{j} \\
  \eta^{j} \\
  \alpha^{j} \\
  \beta^{j} \\
\end{array}%
\right),
\end{equation}
\begin{equation}
\mathbf{W}^{j+1} = \mathbf{T}^{j+1}\mathbf{W}^{j}=
\mathbf{T}^{j+1}\mathbf{T}^{j}\mathbf{W}^{j-1}=\dots=\prod_{k}
\mathbf{T}^{k}\mathbf{W}^{0}.
\end{equation}

Now we will proof a stability with respect to small perturbations of
initial conditions $\mathbf{dW}$ \cite{Halim2002,Halim2003}
\begin{equation}
\mathbf{dW}^{j+1} = \prod_{k} \mathbf{T}^{k}\mathbf{dW}^{0}.
\end{equation}

For the stability conditions we require that the operator $\prod_{k}
\mathbf{T}^{k}$ must be bounded by a constant in a sense of
Frobenius matrix norm
\begin{equation}\label{eq:forbnorm}
\Big|\Big|\prod_{k} \mathbf{T}^{k}\Big|\Big|\leq C.
\end{equation}
For stability, sufficient condition could be write in form
\begin{equation}
||\mathbf{T}^{k}||< \exp(\rho\tau),
\end{equation}
where $\rho$ is a constant, but case with $|\rho(\tau,h)|\leq
const<\infty$
\begin{equation}
||\mathbf{T}^{k}||< \exp(\rho(\tau,h)\tau),
\end{equation}
is also sufficient for stability under condition $\tau, h \to 0$
and some dependence of $\tau$ and $h$.

In matrix $\mathbf{T}$ all elements are matrix with index of
spatial grid. Now we could upper estimate $\rho$. first we upper
estimate all of matrix elements using spectral norm of matrix
\cite{Halim2003}:
\begin{subequations}\label{eq:elemT}
\begin{eqnarray}
  ||T_{11}|| &\leq& 1+\frac{\tau\sigma}{h}, \\
  ||T_{12}|| &\leq& \frac{4\tau k}{h^{2}}+\tau a\max\limits_{i}[(\xi^j_{i})^{2}+(\eta^j_{i})^{2}]+\tau b\max\limits_{i}[(\alpha^j)^{2}_{i}+(\beta^j)^{2}_{i}], \\
  ||T_{21}|| &\leq& \frac{4\tau k}{h^{2}}+\tau a\max\limits_{i}[(\xi^j_{i})^{2}+(\eta^j_{i})^{2}]+\tau b\max\limits_{i}[(\alpha^j)^{2}_{i}+(\beta^j)^{2}_{i}], \\
  ||T_{22}|| &\leq& 1+\frac{\tau\sigma}{h}, \\
  ||T_{33}|| &\leq& 1+\frac{\tau\sigma}{h}, \\
  ||T_{34}|| &\leq& \frac{4\tau k}{h^{2}}+\tau c\max\limits_{i}[(\alpha^j)^{2}_{i}+(\beta^j)^{2}_{i}]+\tau d\max\limits_{i}[(\xi^j)^{2}_{i}+(\eta^j)^{2}_{i}], \\
  ||T_{43}|| &\leq& \frac{4\tau k}{h^{2}}+\tau c\max\limits_{i}[(\alpha^j)^{2}_{i}+(\beta^j)^{2}_{i}]+\tau d\max\limits_{i}[(\xi^j)^{2}_{i}+(\eta^j)^{2}_{i}], \\
  ||T_{44}|| &\leq& 1+\frac{\tau\sigma}{h}.
\end{eqnarray}
\end{subequations}
Taking into account \eqref{eq:forbnorm} and \eqref{eq:elemT} we can
write
\begin{multline}\label{eq:E32}
||\mathbf{T^{j+1}}||\leq 4+4\frac{\tau\sigma}{h}+\frac{16\tau
k}{h^{2}}+\\
2\tau a\max\limits_{i}[(\xi^j_{i})^{2}+(\eta^j_{i})^{2}]+2\tau
b\max\limits_{i}[(\alpha^j)^{2}_{i}+(\beta^j)^{2}_{i}]\\+ 2\tau
c\max\limits_{i}[(\alpha^j)^{2}_{i}+(\beta^j)^{2}_{i}]+2\tau
d\max\limits_{i}[(\xi^j)^{2}_{i}+(\eta^j)^{2}_{i}].
\end{multline}
Taking into account the finite-difference conservation laws
\eqref{E:conserw_law}  we  write
\begin{subequations}
\begin{eqnarray}\label{E:conserw_law1}
I_{u}=\sum_{i=1}^{N}|U_{i}^{j}|^{2}&=&\sum_{i=1}^{N}|U_{i}^{0}|^{2},\\
I_{v}=\sum_{i=1}^{N}|V_{i}^{j}|^{2}&=&\sum_{i=1}^{N}|V_{i}^{0}|^{2},
\end{eqnarray}
\end{subequations}
now we can upper estimate norm \eqref{eq:E32} as
\begin{equation}
||\mathbf{T^{j+1}}||\leq 4+4\frac{\tau\sigma}{h}+\frac{16\tau
k}{h^{2}}+\tau 2(aI_{u}+bI_{v}+dI_{u}+cI_{v}))]\leq
\exp{(\rho\tau)}.
\end{equation}
Finally we can write $\rho$ as
\begin{equation}\label{eq:stability}
\rho = 4\frac{\sigma}{h}+\frac{16 k}{h^{2}}+
4\max(a,d)I_{u}+4\max(b,c)I_{v}.
\end{equation}
This means  that the fraction $\tau/h$ must tend to $0$ (stability
of this system of equation depends on initial energy quantity).For a
better stability this condition requires that $\tau\to 0$ much
faster than $h\to 0$.

\section{Convergence}
In this section we proof that the terms at l.h.s of equations in
finite difference form \eqref{LAB:EXPLICITS} converge to the
equation \eqref{E:CNLSE_2}, which solutions are differentiable with
respect to time and (twice) to $x$. We put into equation
\eqref{LAB:EXPLICITS} a solution in the form
\begin{equation}\label{eq:w}
\mathbf{W}^{j}_{i}=\mathbf{ES}^{j}_{i}+\mathbf{V}^{j}_{i} = \left(%
\begin{array}{c}
  e\xi^{j}_{i}+d \xi^{j}_{i} \\
  e\eta^{j}_{i}+d\eta^{j}_{i} \\
  e\alpha^{j}_{i}+d\alpha^{j}_{i} \\
  e\beta^{j}_{i}+d\beta^{j}_{i}\\
\end{array}%
\right),
\end{equation}
where $\mathbf{ES}$ is an exact solution and $\mathbf{V}$ is a
difference between a numerical solution and the exact solution of
CNLS equations.

Below we show only one of the matrix component
\begin{multline}
\frac{d\xi^{j+1}_{i}-d\xi^{j}_{i}}{\tau}+
\sigma\frac{d\xi^{j}_{i+1}-d\xi^{j}_{i-1}}{2h}+
k\frac{d\eta^{j}_{i+1}-2d\eta^{j}_{i}+d\eta^{j}_{i-1}}{h^{2}}\\
+\frac{e\xi^{j+1}_{i}-e\xi^{j}_{i}}{\tau}+
\sigma\frac{e\xi^{j}_{i+1}-e\xi^{j}_{i-1}}{2h}+
k\frac{e\eta^{j}_{i+1}-2e\eta^{j}_{i}+e\eta^{j}_{i-1}}{h^{2}}\\
+ \left\{a\left[(e\xi^{j}_{i}+d\xi^{j}_{i})^{2}+
(e\eta^{j}_{i}+d\eta^{j}_{i})^{2}\right]+b\left[(e\alpha^{j}_{i}+d\alpha^{j}_{i})^{2}+
(e\beta^{j}_{i}+d\beta^{j}_{i})^{2}\right]\right\}(e\eta^{j}_{i}+d\eta^{j}_{i})=0.
\end{multline}
Let us use the conservation laws \eqref{E:conserw_law} for the
equation \eqref{E:CNLSE_2} and the finite difference equation
\eqref{LAB:EXPLICITS}, that allows us to write
\begin{subequations}
\begin{eqnarray}
I_{eu}&=&\sum_{i=1}^{N}[(e\xi^{j}_{i})^{2}+(e\eta^{j}_{i})^{2}]=\sum_{i=1}^{N}[(e\xi^{0}_{i})^{2}+(e\eta^{0}_{i})^{2}],\\
I_{ev}&=&\sum_{i=1}^{N}[(e\alpha^{j}_{i})^{2}+(e\beta^{j}_{i})^{2}]=\sum_{i=1}^{N}[(e\alpha^{0}_{i})^{2}+(e\beta^{0}_{i})^{2}],
\end{eqnarray}
\end{subequations}
or, in terms of the relation \eqref{eq:w}
\begin{align}
&d\xi^{j+1}_{i}-d\xi^{j}_{i}+
\tau\sigma\frac{d\xi^{j}_{i+1}-d\xi^{j}_{i-1}}{2h}+ \tau
k\frac{d\eta^{j}_{i+1}-2d\eta^{j}_{i}+d\eta^{j}_{i-1}}{h^{2}}\nonumber\\
&+\tau\left\{a\left[(e\xi^{j}_{i})^{2}+
(e\eta^{j}_{i})^{2}\right]+b\left[(e\alpha^{j}_{i})^{2}+
(e\beta^{j}_{i})^{2}\right]\right\}d\eta^{j}_{i}\nonumber
\\
&-\tau\left\{a\left[(e\xi^{j}_{i})^{2}+
(e\eta^{j}_{i})^{2}\right)+b\left[(e\alpha^{j}_{i})^{2}+
(e\beta^{j}_{i})^{2}\right]\right\}d\eta^{j}_{i}\nonumber
\\%
&-\tau\left\{a\left[(e\xi^{j}_{i})^{2}+
(e\eta^{j}_{i})^{2}\right]+b\left[(e\alpha^{j}_{i})^{2}+
(e\beta^{j}_{i})^{2}\right]\right\}e\eta^{j}_{i}\nonumber
\\
&+\tau\left\{a\left[(e\xi^{j}_{i}+d\xi^{j}_{i})^{2}+
(e\eta^{j}_{i}+d\eta^{j}_{i})^{2}\right]+b\left[(e\alpha^{j}_{i}+d\alpha^{j}_{i})^{2}+
(e\beta^{j}_{i}+d\beta^{j}_{i})^{2}\right]\right\}(e\eta^{j}_{i}+d\eta^{j}_{i})=\nonumber\\
&-\tau\left(\frac{e\xi^{j+1}_{i}-e\xi^{j}_{i}}{\tau}-
\sigma\frac{e\xi^{j}_{i+1}-e\xi^{j}_{i-1}}{2h}-
k\frac{e\eta^{j}_{i+1}-2e\eta^{j}_{i}+e\eta^{j}_{i-1}}{h^{2}}\right)\nonumber\\
&+\tau\left\{a\left[(e\xi^{j}_{i})^{2}+
(e\eta^{j}_{i})^{2}\right]+b\left[(e\alpha^{j}_{i})^{2}+
(e\beta^{j}_{i})^{2}\right]\right\}e\eta^{j}_{i}.
\end{align}
The right side of the equation for the differentiable solutions is
of the order $\Theta(\tau + h+h^{2})$
\cite{book:Hoffman:numerical_eng}, while to represent the left side
we use the matrix $\mathbf{T}$. Hence one arrives at
\begin{multline}
d\xi^{j+1}_{i}=\mathbf{T}^{j+1}_{1i}\mathbf{V}_{i}^{j}
+\tau\left\{a\left[(e\xi^{j}_{i})^{2}+
(e\eta^{j}_{i})^{2}\right]+b\left[(e\alpha^{j}_{i})^{2}+
(e\beta^{j}_{i})^{2}\right]\right\}(e\eta^{j}_{i}+d\eta^{j}_{i})
\\%
-\tau\left\{a\left[(e\xi^{j}_{i}+d\xi^{j}_{i})^{2}+
(e\eta^{j}_{i}+d\eta^{j}_{i})^{2}\right]+b\left[(e\alpha^{j}_{i}+d\alpha^{j}_{i})^{2}+
(e\beta^{j}_{i}+d\beta^{j}_{i})^{2}\right]\right\}(e\eta^{j}_{i}+d\eta^{j}_{i})\\
+\Theta(\tau + h+h^{2}).
\end{multline}
Let us define the norm of the numerical solution as
\begin{equation}
\|U^{j}\|=h\left(\sum_{i}|U_{i}^{j}|^{2}\right)^{\frac{1}{2}}.
\end{equation}

Now we  upper estimate one element of the vector
 $\mathbf{V}^{j}$
\begin{align}
\|d\xi^{j+1}\|&\leq\|\mathbf{T}^{j+1}_{1}\|\|\mathbf{V}^{j}\|
+\tau\left(aI_{u}+bI_{v}\right)\|\eta^{j}\|
-\tau\left(aI_{ue}+bI_{ve}\right)\|\eta^{j}\|+ \Theta(\tau + h+h^{2})\\
&\leq\exp{(j\tau\rho)}\|\mathbf{V}^{0}\|+\tau\left(aI_{u}+bI_{v}\right)\|U^{j}\|
-\tau\left(aI_{ue}+bI_{ve}\right)\|U^{j}\|+\Theta(\tau + h+h^{2}),
\end{align}

Next we  write the converge condition for all $\mathbf{V}^{j+1}$
matrix components using the Frobenius norm up to the choice of the
initial error  $||\mathbf{V}^{0}||=0$.

\begin{align}
\|\mathbf{V}^{j+1}\|&\leq Q+\Theta(\tau +h+h^{2}),
\end{align}
where $i=1,2,3,4$ and we define
\begin{align}
Q =4|
\max{(a,b,c,d)}(I_u+I_{v})^{3/2}-\max{(a,b,c,d)}(I_{ue}+I_{ve})^{3/2}|.
\end{align}

\section{Numerical results}

\subsection{Nonlinear Sch\"{o}dinger equation}
\begin{figure}[htp]
  \centering
  \subfigure[]{\includegraphics[width=0.45\textwidth]{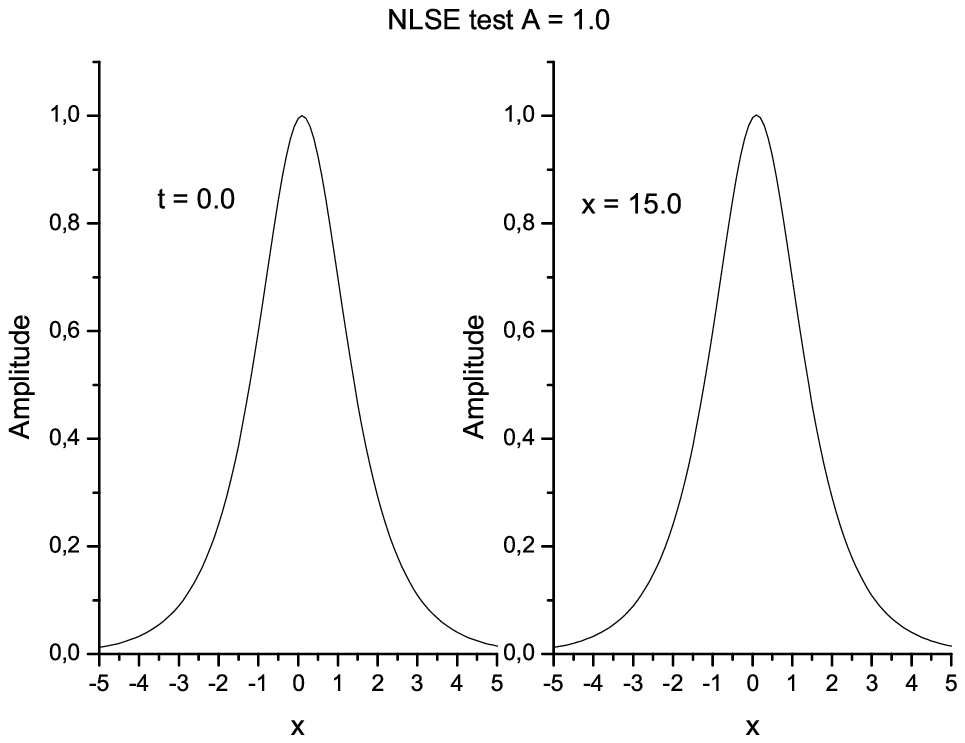}}
  \subfigure[]{\includegraphics[width=0.45\textwidth]{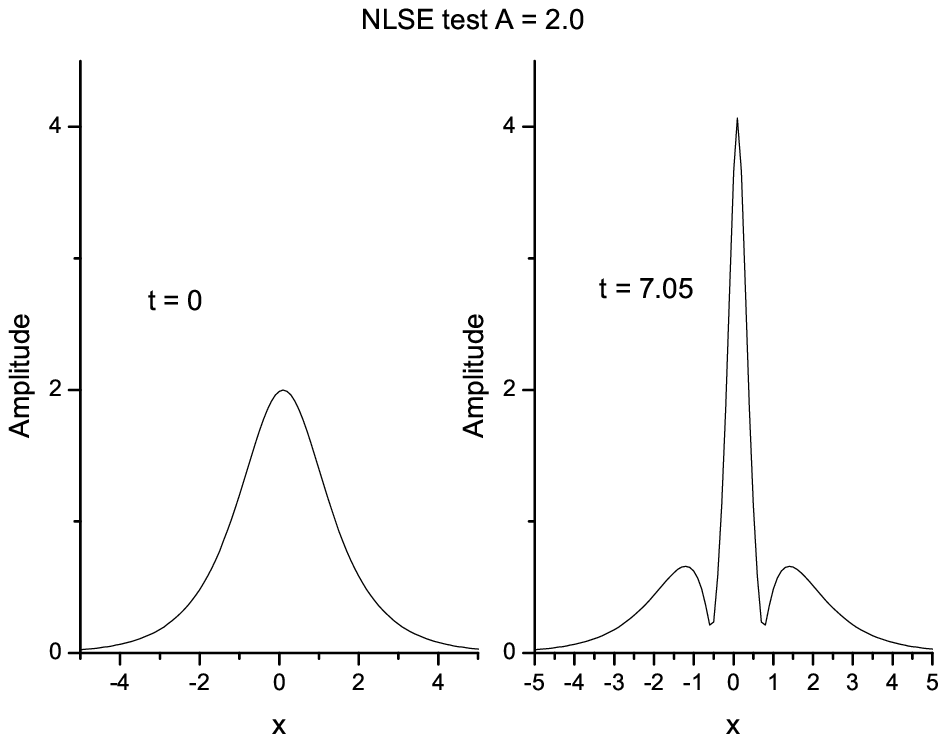}}
  \caption{NLS equation with initial amplitude ($A_1$) a) 1 and b) 2}
  \label{picLAB:NLS1}
\end{figure}

For simple test we put $c=d=0$ and $k=0.5$ in this case we have
simple nonlinear Schr\"{o}dinger equation. First we test NLS
equation with initial condition $U(0,x)=\mathrm{sech}(x)$ which
pulse should  not change shape during propagation (picture
\ref{picLAB:NLS1}a).

As a full initial condition we set:
\begin{subequations}
\begin{eqnarray}
U(0,x)&=&A_1 \mathrm{sech}(x),\\
V(0,x)&=&0,
\end{eqnarray}
\end{subequations}
this lead to solve NLSE problem in form
\begin{equation}\label{LAB:NLSEq}
iU_{t}+\frac{1}{2}U_{xx}+|U|^{2}U=0.
\end{equation}

When amplitude $A_1 = 1$ the exact solution of equation
\eqref{LAB:NLSEq} is
\begin{equation}
U(t,x)=A_1 \mathrm{sech}(x)\exp{(it/2)},
\end{equation}

\begin{figure}[htp]
  \centering
  \subfigure[]{  \includegraphics[width=0.45\textwidth]{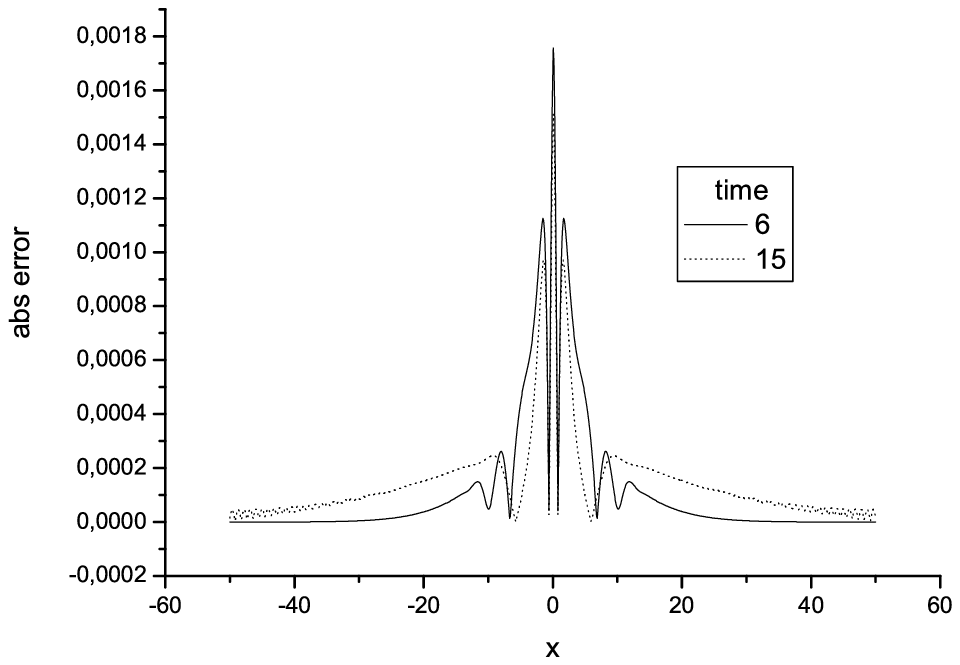}}
  \subfigure[]{  \includegraphics[width=0.45\textwidth]{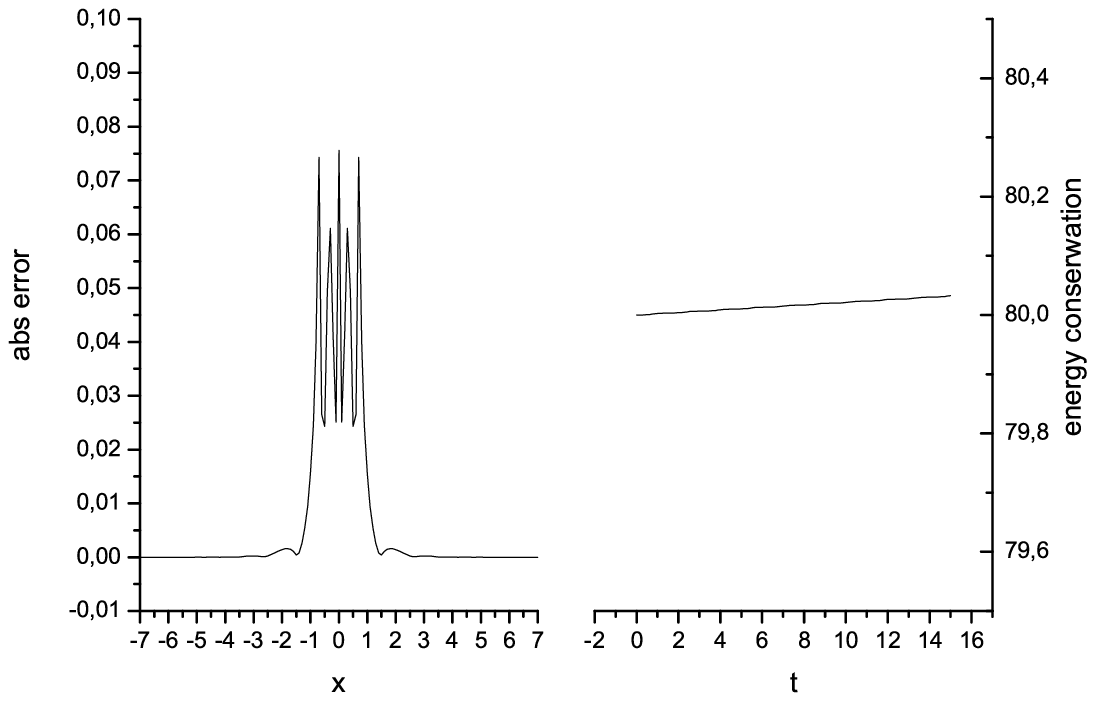}}
  \caption{Error of numerical calculation for different time value a) $A_{1}=1$ and b) $A_{1}=2$}
  \label{LAB:NLS1_ER}
\end{figure}

Next we test NLS equation with amplitude $A_1 = 2$; the results are
the same as in \cite{book:agrawal:NonFibOpt} (picture
\ref{picLAB:NLS1}b). For this case the exact solution of the
equation \eqref{LAB:NLSEq} takes more complicated form
\cite{book:agrawal:NonFibOpt}
\begin{equation}
U(t,x)=\frac{[\cosh(3x)+\exp(4it)\cosh(x)]\exp(it/2)}{\cosh(4x)+4\cosh(2x)+3\cos(4t)}.
\end{equation}
The numerical solution for the second case deviates more than the
first one because it has bigger energy per pulse  and we use for
both cases the same time and space steps (see equation
\ref{eq:stability}).

\subsection{Manakov solitons}

\begin{figure}
\begin{minipage}[t]{7.5cm}
\begin{center}
\includegraphics[width=7.5cm,clip]{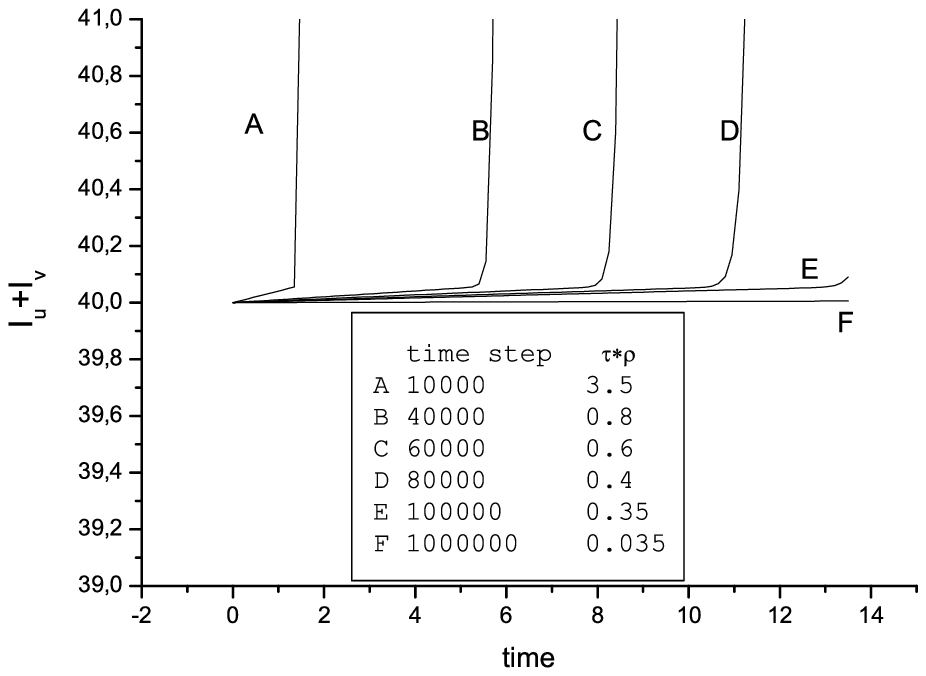}
\caption{\label{LAB:STABMAN} Stability of explicit method for
Manakov examples}
\end{center}
\end{minipage}
\hfill\hfill
\begin{minipage}[t]{7.5cm}
\begin{center}
\includegraphics[width=7.5cm,clip]{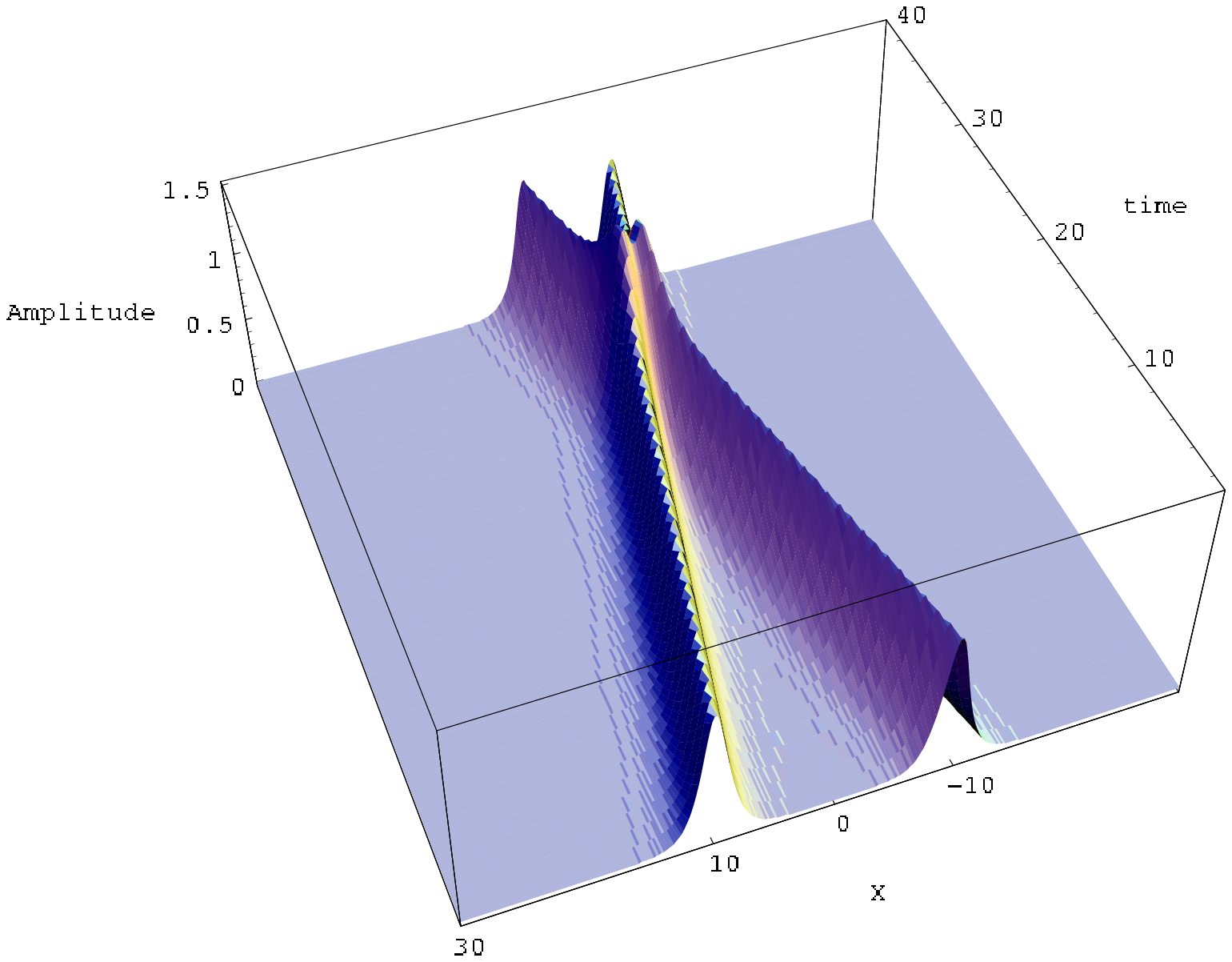}
\caption{\label{LAB:ART1} Manakov case:$A_1=A_2 = 1$, nonlinear
coefficients ($a$,$b$,$c$,$d$) = 1} and $\sigma=1$.
\end{center}
\end{minipage}
\end{figure}

Let us consider the soliton solutions of the Manakov system as
examples to test the stability and the convergence of the explicit
scheme. We start with the stability of the explicit scheme. At the
picture \ref{LAB:STABMAN}   six cases of energy conservation
behaviour  are shown for the values of $\rho\tau$  which are bigger
than $0.1$: we observe very unstable results. When we decrease the
time step, our solution is stabilized. It is very important to
remark that $\rho$ depends on initial amplitude of pulses. We choose
the initial conditions as
\begin{subequations}
\begin{eqnarray}
U(0,x)&=&A_1 \mathrm{sech}(x),\\
V(0,x)&=&A_2 \mathrm{sech}(x).
\end{eqnarray}
\end{subequations}

\begin{table}[ht]
\caption{Parameters for the numerical experiment ''Manakov
solitons''} \centering
\begin{tabular}{|l|l|}
\hline
x & -50 ... 50 \\
\hline
time & 15 \\
\hline
sigma & 0 \\
\hline
a & 1 \\
\hline
b & 1 \\
\hline
c & 1 \\
\hline
d & 1 \\
\hline
$A_1$ & 1 \\
\hline
$A_2$ & 1 \\
\hline
space steps & 1000 \\
\hline
time steps & 10000 - 1000000 \\
\hline
\end{tabular}
\end{table}

\subsection{Collision of two solitons}
We make this experiment to compare the results with the ''experiment
1'' of the paper \cite{Sun2004}. General properties of such
collisions of two solitons are well known \cite{Porsezian2001},
hence we do not focus on details of this type of soliton
interaction.

\begin{figure}[htp]  \label{LAB:ART2}
  \centering
  \subfigure[$|U|$]{  \includegraphics[width=0.45\textwidth]{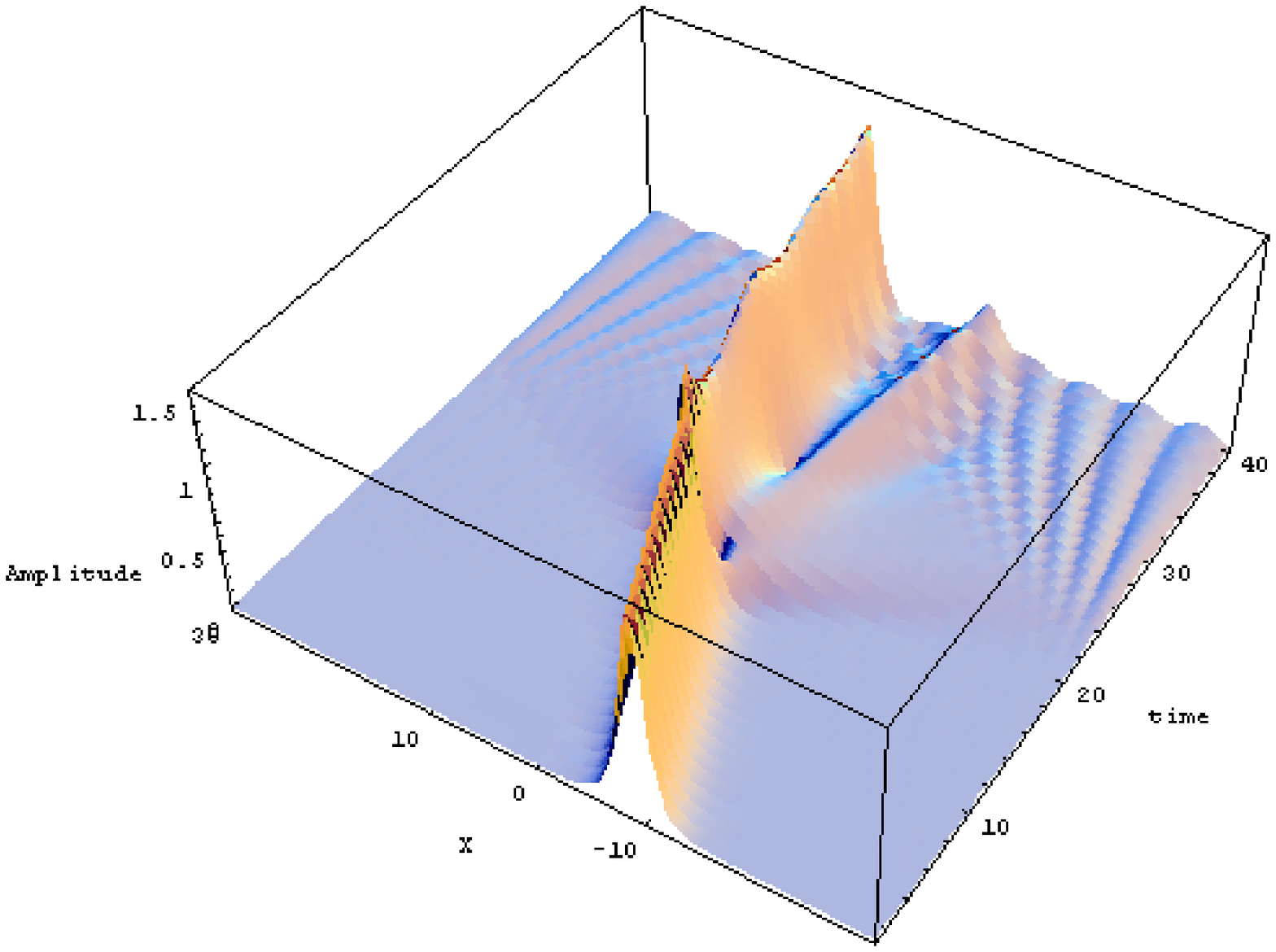}}
  \subfigure[$|V|$]{  \includegraphics[width=0.45\textwidth]{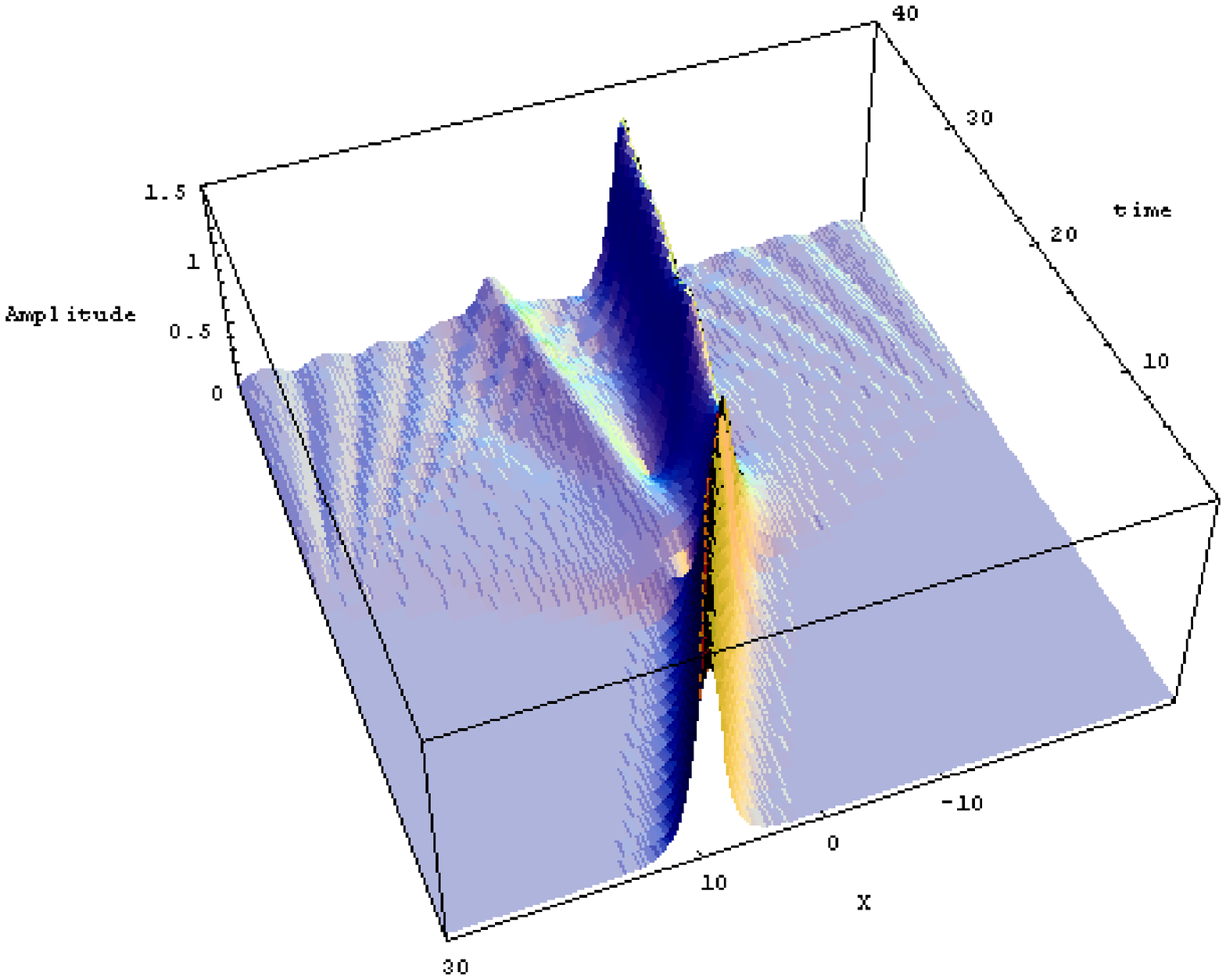}}
  \caption{Inelastic collision of two solitons}
\end{figure}

\subsection{Different group velocities}

The examples to be studied in this section are most important for
us, because we interest in describing modes interaction in a
nonlinear waveguide of different modes excited (different modes have
different group velocitys) \cite{Snyder1978,Leble2005}. This
investigation could be useful not only for the situation described
below, but for soliton trains interactions as well
\cite{Mamyshev1991}.

 As initial conditions
we took two $\mathrm{sech}$-impulses and each equation have
different nonlinear coefficients. This could appear in a waveguide
when two modes are excited with different group velocity. We set
initial conditions for this case as:
\begin{subequations}
\begin{eqnarray}
U(0,x)=A_1 \mathrm{sech}(x+D_1)\exp{(iv_1 x)},\\
V(0,x)=A_2 \mathrm{sech}(x+D_2)\exp{(iv_2 x)}.
\end{eqnarray}
\end{subequations}

\begin{figure}[htp]
  \centering
  \subfigure[]{\includegraphics[width=0.45\textwidth]{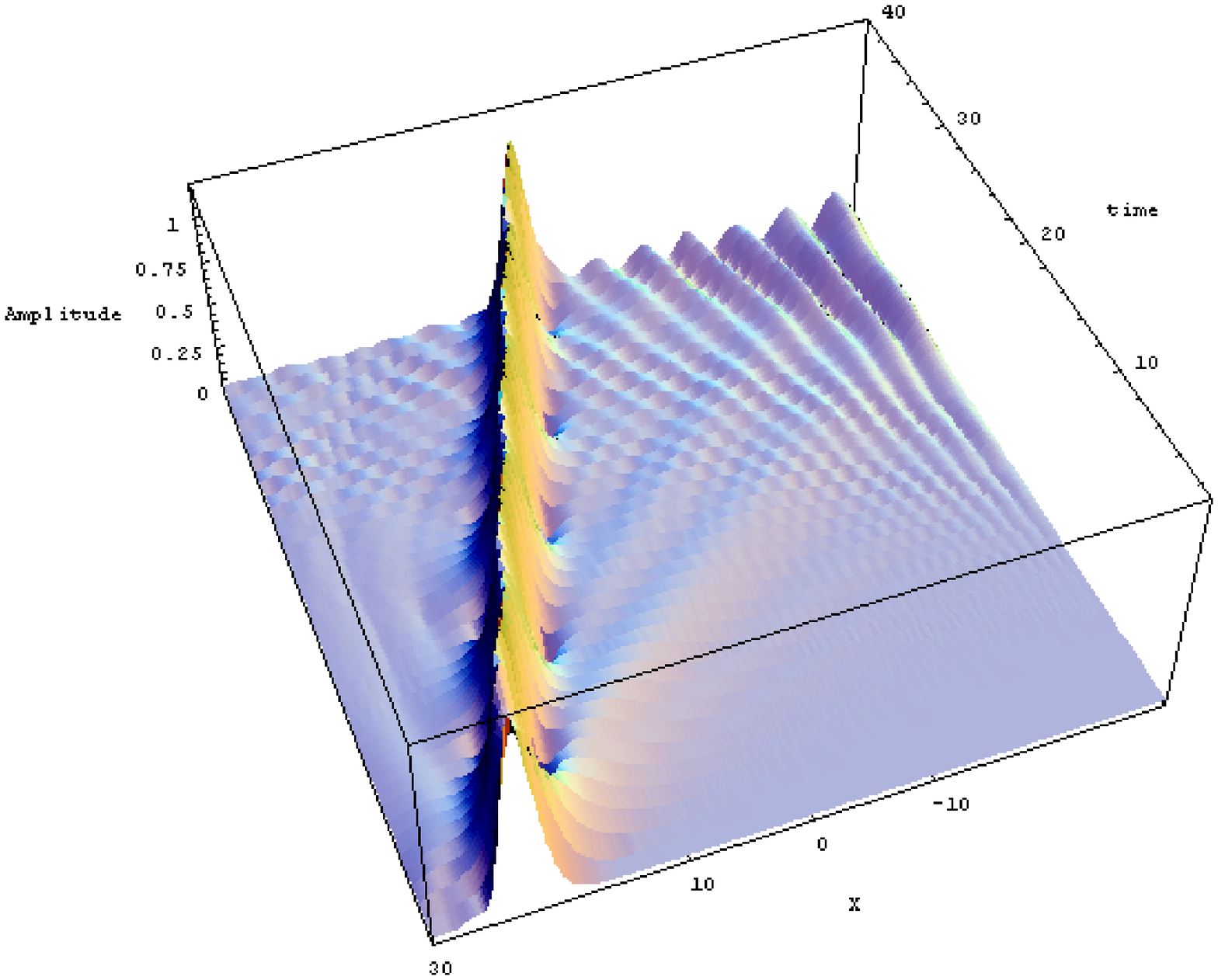}}
  \subfigure[]{\includegraphics[width=0.45\textwidth]{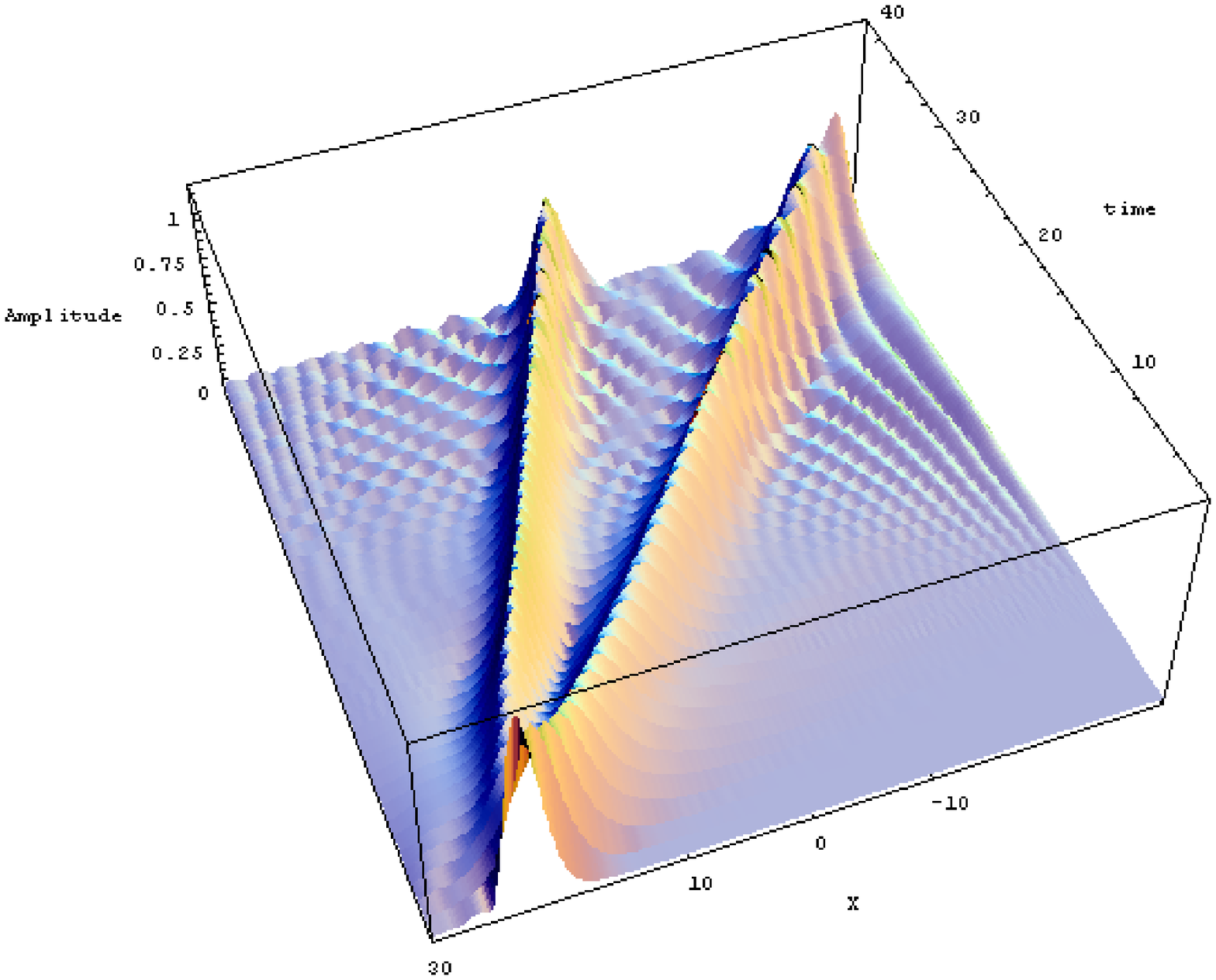}}
  \caption{Two cases of impulses with different group velocity (for parameters see table \ref{TAB:DIFF})}
  \label{LAB:DIF12}
\end{figure}

Consider two solitons one with zero velocity and second with
velocity greater than zero, we show this situations on picture
\ref{LAB:DIF12}. This two cases of pulses interaction differ only
by group velocity (parameters for this pulses are in the table
\ref{TAB:DIFF}).

 This two impulses start from the same
position (e.g. two modes excited in the waveguide). When velocity
of second pulse is $0.7$ (for given parameters) first impulses
intercept part of energy thought nonlinear interactions and move
with second pulse with average velocity for both pulses (picture
\ref{LAB:DIF12}a). When the velocity of second impulses is high
enough we have two pulses which moves with different group
velocity see picture \ref{LAB:DIF12}b (nonlinear interaction
between pulses happen only at the start of propagation).

\begin{table}[ht] \label{TAB:DIFF}
\caption{Parameters for numerical experiment ''Different group
velocity''} \centering
\begin{tabular}{|l|l|l|}
\hline
& (a)&(b)\\
 \hline
x & -30 ... 30& \\
\hline
time & 40 &40\\
\hline
$\sigma$ & 0&0 \\
\hline
a & 1& 1 \\
\hline
b & 1/3& 1/3 \\
\hline
c & 1& 1 \\
\hline
d & 1/3& 1/3 \\
\hline
$A_1$ & 1.2& 1.2 \\
\hline
$A_2$ & 1.4 & 1.4 \\
\hline
$v_1$ & 0.7& 0.95\\
\hline
$v_2$ & 0& 0\\
\hline
h & 0.2& 0.2 \\
\hline
$\tau$ & 0.02& 0.02 \\
\hline
\end{tabular}
\end{table}

\subsection{Difference between explicit and implicit scheme}
\begin{figure}[htp]
  \centering
  \subfigure[]{\includegraphics[width=0.45\textwidth]{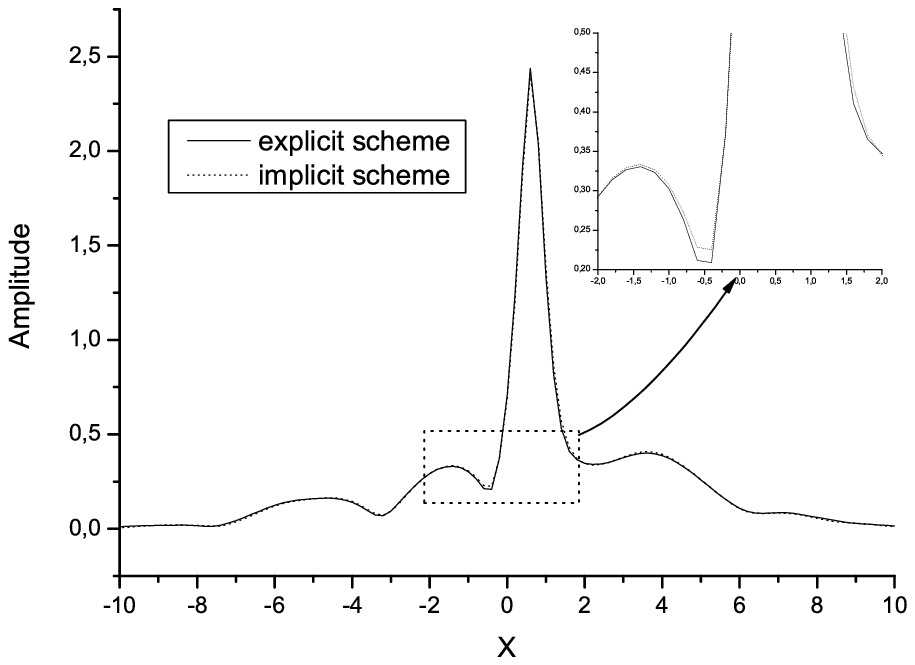}}
  \subfigure[]{\includegraphics[width=0.45\textwidth]{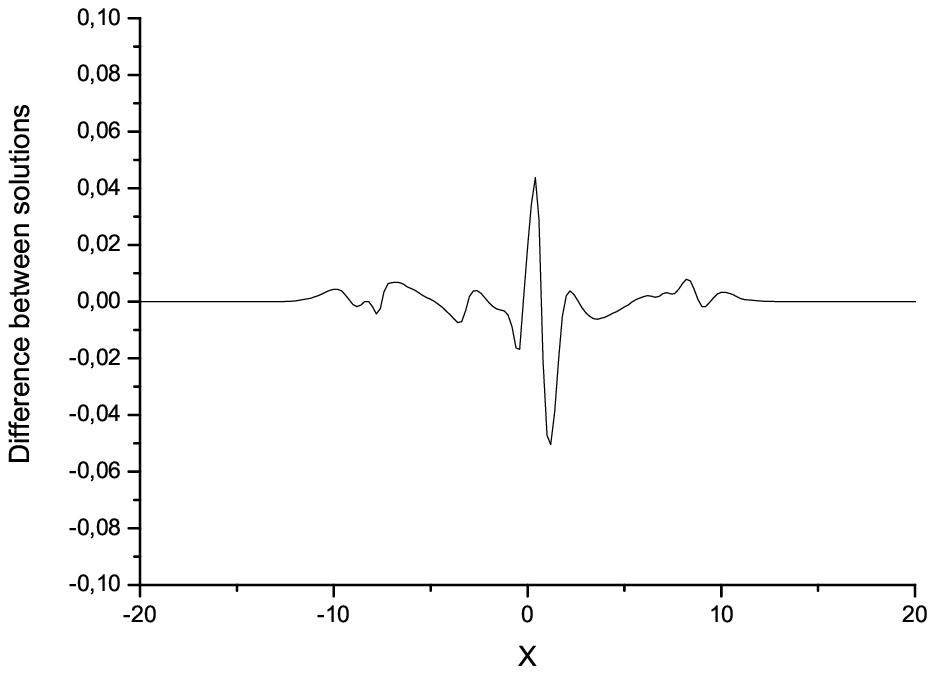}}
  \caption{Compare between explicit and implicit method.}
  \label{LAB:EXPIMPL}
\end{figure}

If we would like to have more space details (smaller $h$), in
explicit scheme we must take adequate smaller time step to assure
stability of numerical scheme. This have very big influence on
calculation time. In an implicit scheme we could use an iteration
method \cite{book:Hoffman:numerical_eng,Chang1999} and  shorter time
of calculation is achieved. On the picture \ref{LAB:EXPIMPL} we show
numerical results for both methods. Note, the time step for the
explicit scheme is much larger than for the implicit scheme but each
step in the implicit scheme needs 2-4 iterations.

\begin{table}[ht]
\caption{Parameters for numerical experiment for implicit and
explicit method} \centering
\begin{tabular}{|l|l|l|}
\hline
Parameter & explicit & implicit \\
\hline
space step & 300 & 300 \\
\hline
time step & 1000000 & 2000 \\
\hline
time & 40 & 40 \\
\hline
x & -30...30 & -30...30 \\
\hline
$A_1$ & 1.5 & 1.5 \\
\hline
$A_2$ & 1.5 & 1.5 \\
\hline
a & 1 & 1 \\
\hline
b & 0.2 & 0.2 \\
\hline
c & 1 & 1 \\
\hline
d & 1.6 & 1.6 \\
\hline
$\sigma$ & 0.3 & 0.3 \\
\hline
\end{tabular}
\end{table}

\subsection{Rectangular pulse decay}
\begin{figure}[htp]
  \centering
  \subfigure[]{\includegraphics[width=0.45\textwidth]{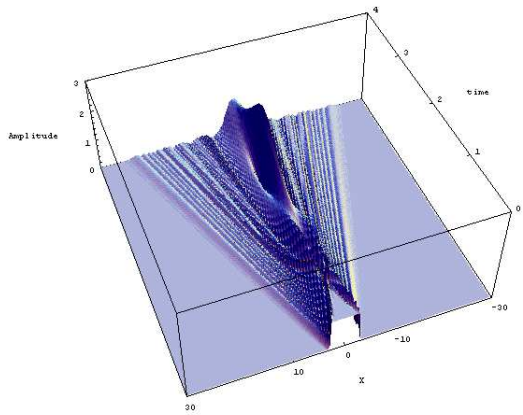}}
  \subfigure[]{\includegraphics[width=0.45\textwidth]{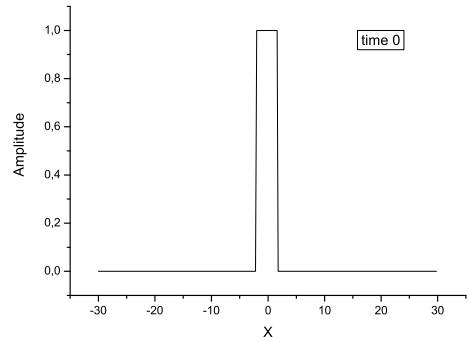}}\\
  \subfigure[]{\includegraphics[width=0.45\textwidth]{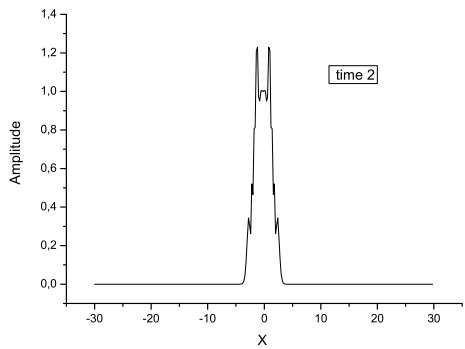}}
  \subfigure[]{\includegraphics[width=0.45\textwidth]{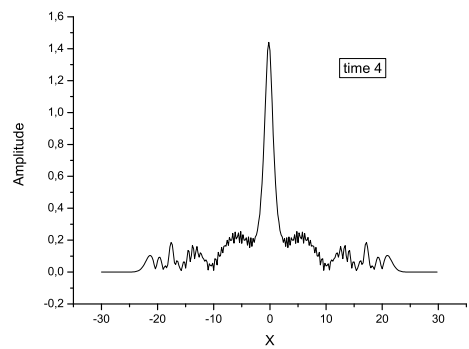}}
  \caption{Rectangular pulse evolution. On picture a) there is a 3D plot on picture b),c) and d) there are intersection of plot a) in different time. }
  \label{picLAB:PROST}
\end{figure}

In this subsection we engaged in asymptotic solution of  Nonlinear
Schr\"{o}dinger Equation. As initial condition we choose rectangular
pulse as in \cite{book:novikov:solitontheory} (see picture
\ref{picLAB:PROST}b).

\begin{figure}[htp]
  \centering
  \subfigure[]{\includegraphics[width=0.45\textwidth]{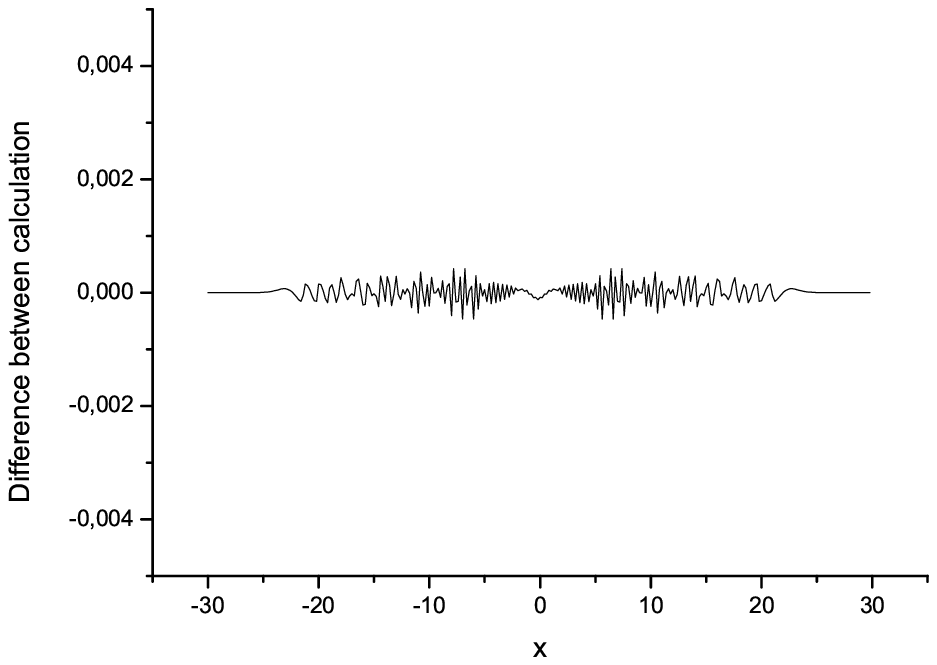}}
  \subfigure[]{\includegraphics[width=0.45\textwidth]{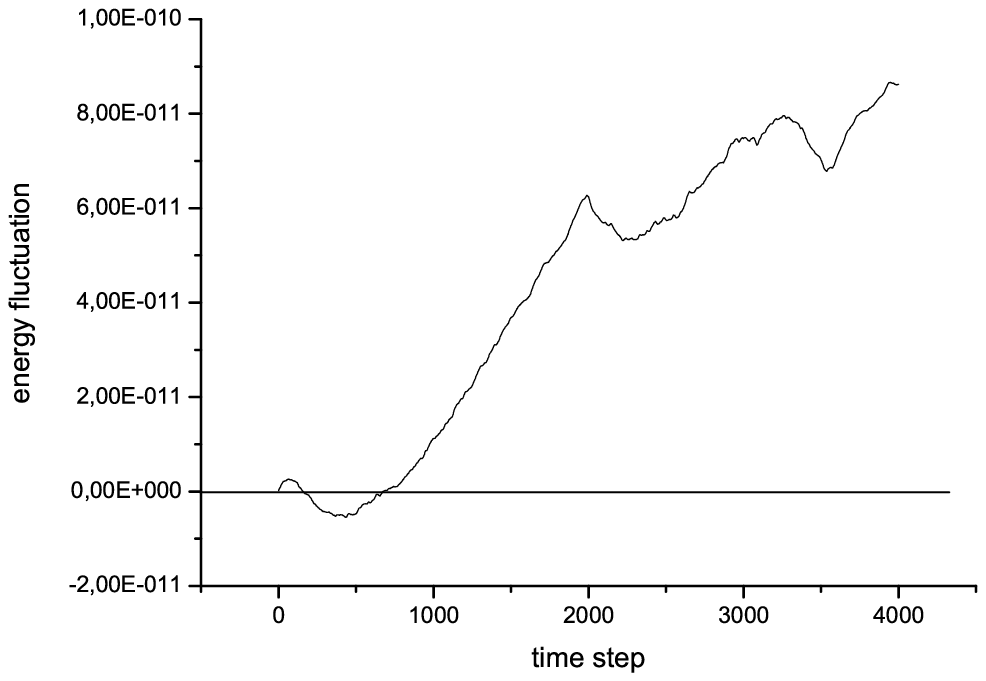}}
  \caption{\textbf{a)} Difference between two numerical solutions
   with different time steps ($\tau=2e^{-4}$ and $\tau=2e^{-5}$.
  \textbf{b)} The energy deviation (to the initial state) for the time step $\tau=2e^{-5}$.}
  \label{picLAB:ENEPRO}
\end{figure}

We test our method for this kind of initial pulse. First we made two
calculations with different time steps (for the first $\tau=2e^{-4}$
and $\tau=2e^{-5}$ for the second one) and compare errors between
(picture \ref{picLAB:ENEPRO}a). The difference between calculations
is of the order determined by $\tau\rho$. Next we check the energy
conservation (see picture \ref{picLAB:ENEPRO}b).

\begin{figure}[htp]
  \centering
  \includegraphics[width=0.75\textwidth]{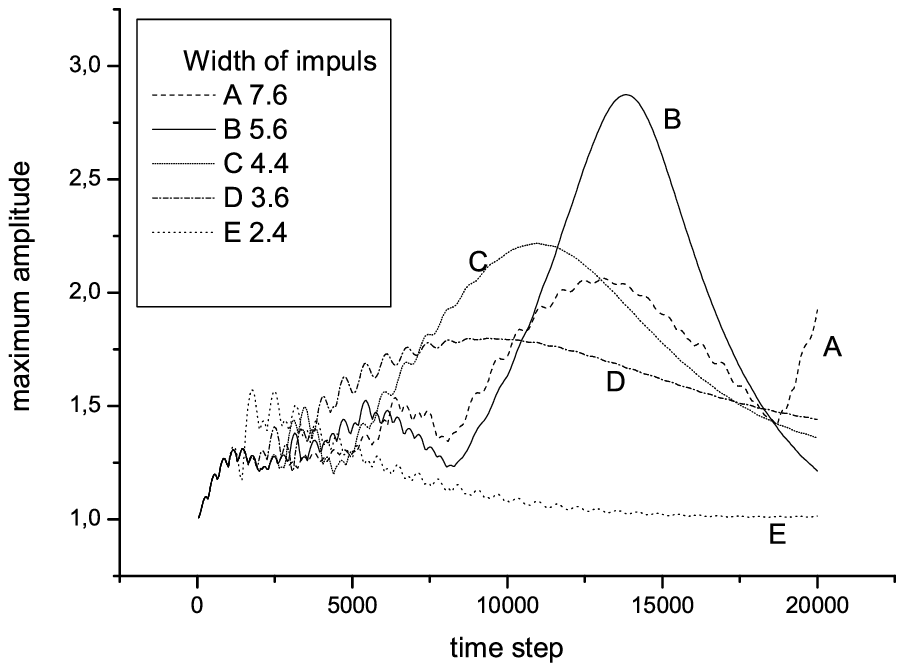}
  \caption{Maximum amplitude for different pulse width (in round brackets) for time step $\tau=2e^{-4}$ (maximum time in dimensionless unit is equal 4).}
  \label{picLAB:IMPW}
\end{figure}

For the presented time ($t=4$) we expect that the transient stage
from rectangular shape of the pulse to the asymptotic behaviour of
solution is realized. On the picture \ref{picLAB:IMPW} the evolution
(maximum amplitude) of rectangular pulses with different width is
presented. If the pulse have smaller width then the decay is faster.

\section{Conclusion}
We study a convergence and stability conditions which would be
useful to estimate other numerical schemes for NLS as well.
Additionally these conditions could be compared with ones for
numerical methods for other nonlinear equations \cite{Halim2003} and
show a character of CNLS equations (from numerical side). Most
important results which we show in this paper is how to estimate a
time step for the explicit method due to initial energy of pulses.
The results could be adapted to other methods based on this basic
scheme which we use \cite{Ismail2004}. There is a possibility to
enhance this method to third or high order in time and this method
could be used as a starting point (to calculate lacked points of
grid).

\end{document}